\documentclass[12pt]{article}
\usepackage{amsfonts}
\usepackage{mathrsfs}
\usepackage{amsmath,amssymb}
\baselineskip 5.5pt \lineskip 5.5pt \numberwithin{equation}{section}

\def\M{\mathcal {M}}
\def\V{\mathcal {V}}

\def\l{\lambda}

\def\LL{\mathscr{L}}

\def\cl{\centerline}

\def\vs{\vspace*}

\def\C{\mathbb{C}}

\def\Z{\mathbb{Z}}

\newtheorem{theo}{Theorem}[section]
\newtheorem{conj}[theo]{Conjecture}
\newtheorem{prob}[theo]{Problem}

\newtheorem{coro}[theo]{Corollary}

\newtheorem{lemm}[theo]{Lemma}

\begin{document}
\cl{\large\bf Classification of Irreducible Weight Modules with a}

\cl{\large\bf Finite-dimensional Weight Space over the}

\cl{\large\bf Twisted Schr\"{o}dinger-Virasoro Lie algebra}\vs{8pt}

\cl{{\bf Junbo Li}}\vs{4pt} \cl{\small\it Department of Mathematics,
Shanghai Jiao Tong University Shanghai 200240, P. R. China;}

\cl{\small\it Department of Mathematics, Changshu Institute of
Technology, Changshu 215500, P. R.
China}
\cl{\small\it Email: sd$\_$junbo@163.com}\vs{10pt} \cl{{\bf Yucai
Su}\footnote{Supported by NSF grants 10471091 of China and ``One
Hundred Talents Program'' from University of Science and Technology
of China}} \cl{\small\it Department of Mathematics, University of
Science and Technology of China} \cl{\small\it Hefei 230026, P. R.
China } \cl{\small\it Email: ycsu@ustc.edu.cn}\vs{10pt}

\noindent{\small {\bf Abstract.} It is shown that the support of an
irreducible weight module over the Schr\"{o}dinger-Virasoro Lie
algebra with an infinite-dimensional weight space, coincides with
the weight lattice and that all nontrivial weight spaces of such a
module are infinite-dimensional. As a side-product, it is obtained
that every simple weight module over the Schr\"{o}dinger-Virasoro
Lie algebra with a nontrivial finite-dimensional weight space, is
a Harish-Chandra module.\\[7pt]
\noindent{\bf Key Words:} The Schr\"{o}dinger-Virasoro Lie algebra,
weight modules, irreducible modules\\
{\it Mathematics Subject Classification (2000)}: 17B10, 17B65,
17B68.} \vs{10pt}\par \cl{\bf1. \ Introduction}
\setcounter{section}{1}\setcounter{theo}{0} The {\it twisted
Schr\"{o}dinger-Virasoro Lie algebra} is the infinite-dimensional
Lie algebra $\LL$ with  $\C$-basis $\{L_n,Y_m,M_p\,|\,m,n,p\in \Z\}$
and the following relations
\begin{eqnarray}
\!\!\!&\!\!\!& [L_n,L_p]=(p-n)L_{n+p},\nonumber
\\[4pt]\!\!\!&\!\!\!&
[L_n,Y_m]=(m-\frac{n}{2})Y_{n+m},\nonumber\\[4pt]\!\!\!&\!\!\!&
[L_n,M_p]=pM_{n+p},
\nonumber\\[4pt] \!\!\!&\!\!\!&
[Y_m,Y_{m'}]=(m'-m)M_{m+m'},
\nonumber\\[4pt]\!\!\!&\!\!\!&
[Y_m,M_p]=[M_n,M_p]=0.\nonumber
\end{eqnarray}
where $m,m',n,p\in\Z$. This infinite-dimensional Lie algebra is the
twisted deformation of the Schr\"{o}dinger-Virasoro Lie algebra,
which was originally introduced in [\ref{H1}] in the context of
non-equilibrium statistical physics, as a by-product of the
computation of $n$-point functions that are covariant under the
action of the Schr\"{o}dinger group, containing as subalgebras both
the Lie algebra of invariance of the free Schr\"{o}dinger equation
and the central charge-free Virasoro algebra. Both original and
twisted Schr\"{o}dinger-Virasoro Lie algebras are closely related to
the Schr\"{o}dinger Lie algebras and the Virasoro Lie algebra (cf.
[\ref{DZ},\ref{HU}-\ref{M},\ref{SS},\ref{S},\ref{Z}]). They should
consequently play a role akin to that of the Virasoro Lie algebra in
two-dimensional equilibrium statistical physics.

One attempt of introducing two-dimensional conformal field theory is
to understand the universal behaviour of two-dimensional statistical
systems at equilibrium and at the critical temperature, where they
undergo a second-order phase transition. A systematic investigation
of the theory of representations of the Virasoro algebra in the
80¡¯es led to introduce a class of physical models (called unitary
minimal models), corresponding to the unitary highest weight
representations of the Virasoro algebra with central charge less
than one. Miraculously, covariance alone is enough to allow the
computation of the statistic correlators (or so-called `$n$-point
functions') for these highly constrained models.

Motivated by the search for deformations and central extensions of
both original and twisted Schr\"{o}dinger-Virasoro Lie algebras,
Roger and Unterberger presented the sets of generators provided by
the cohomology classes of the cocycles in the fifth part of
[\ref{RU}]. In [\ref{U}], Unterberger constructed vertex algebra
representations of the Schr\"{o}dinger-Virasoro Lie algebras out of
a charged symplectic boson and a free boson and its associated
vertex operators.

A {\it weight module} $\M$ over ${\cal L}$ is an $\frak
H$-diagonalizable module, where $${\frak
H}:=\mbox{Span}_{\C}\{L_0,M_0\}$$
 is the {\it Cartan subalgebra} of $\LL$. If, in
addition, $\M$ is irreducible and all weight spaces $\M_\l$ (cf.
(\ref{eq1.2})) are finite-dimensional, the module is called a {\it
Harish-Chandra module}.

If $\M$ is an irreducible weight $\LL$-module, then $M_0$ must act
as some complex number $h_M$ on $\M$. Furthermore, $\M$ has the
weight space decomposition
\begin{eqnarray}
\label{eq1.2} \M=\bigoplus_{\lambda\in\C}\M_\lambda, \mbox{ \ \
where } \M_\lambda=\{v\in \M\,|\,L_0v=\lambda v\},
\end{eqnarray}
where $\M_\l$ is called {\it a weight space} with weight $\l$.
Denote the set of {\it weights $\l$} of $\M$ by
$${\rm Supp}(\M):=\{\lambda\in\C\,|\,\M_\lambda\neq 0\},$$
called the {\it support} of $\M$. Obviously, if $\M$ is an
irreducible weight $\LL$-module, then there exists $\lambda\in\C$
such that ${\rm Supp}(\M)\subset\lambda+\Z$.

An irreducible weight module $\V$ is called {\it a pointed module}
if there exists a weight $\l\in\C$ such that ${\rm dim\,}\V_\l=1$.
Xu posted the following natural problem in [\ref{X}]:
\begin{prob}
\label{prob1.1} Is any irreducible pointed module over the Virasoro
algebra a Harish-chandra module?
\end{prob}

An irreducible weight module $\V$ is called {\it a  mixed module} if
there exist $\l\in\C$ and $i\in\Z$ such that ${\rm
dim\,}\V_\l=\infty$ and ${\rm dim\,}\V_{\l+i}<\infty$. The following
conjecture was posted in [\ref{M}]:
\begin{conj}
\label{conj1.2} There are no irreducible mixed modules over the
Virasoro algebra.
\end{conj}

Mazorchuk and Zhao [\ref{MZ}] gave the positive answers to the above
question and conjecture. In [\ref{SS}], the authors gave positive
answers to the above question and conjecture for the twisted
Heisenberg-Virasoro algebra. In this note, we also give the positive
answers to the above question and conjecture for the
Schr\"{o}dinger-Virasoro Lie algebra. Our main result is the
following:
\begin{theo}
\label{theo1.3} Let $\M$ be an irreducible weight $\LL$-module.
Assume that there exists $\lambda\in\C$ such that ${\rm
dim\,}\M_\lambda=\infty$. Then ${\rm Supp}(\M)=\lambda+\Z$ and for
every $k\in\Z$, we have ${\rm dim\,}\M_{\lambda+k}=\infty$.
\end{theo}

Theorem \ref{theo1.3} also implies the following classification of
all irreducible weight $\LL$-modules which admit a nontrivial finite
dimensional weight space:
\begin{coro}
\label{coro1.4} Let $\M$ be an irreducible weight $\LL$-module.
Assume that there exists $\lambda\in\C$ such that $0<{\rm
dim\,}\M_\lambda<\infty$. Then $\M$ is a Harish-Chandra module.
\end{coro}
\cl{\bf2. \ Proof of Theorem \ref{theo1.3}}
\setcounter{section}{2}\setcounter{theo}{0} \setcounter{equation}{0}
\vs{1pt}\par

\noindent We first recall a main result about the weight
Virasoro-modules in [\ref{MZ}].
\begin{theo}
\label{theo2.1} Let $\V$ be an irreducible weight Virasoro-module.
Assume that there exists $\lambda\in \C$, such that ${\rm
dim\,}\V_\lambda=\infty$. Then ${\rm Supp}(\V)=\lambda+\Z$ and for
every $k\in\Z$, we have ${\rm dim\,}\V_{\lambda+k}=\infty$.
\end{theo}
\begin{lemm}
\label{lemm2.2} Assume that there exists $\mu\in\C$ and a non-zero
element $v\in\M_\mu$, such that
\begin{eqnarray*}
&&L_1v=L_2v=Y_{1}v=M_1v=0\ \ \ \mbox{or}\\
&&L_{-1}v=L_{-2}v=Y_{-1}v=M_{-1}v=0.
\end{eqnarray*}
Then $\M$ is a Harish-Chandra module.
\end{lemm}
{\it Proof.~}~Indeed, under these conditions, $v$ is either a
highest or a lowest weight vector and hence $\M$ is a Harish-Chandra
module (see, e.g. [\ref{LZ}]).$\hfill\Box$

Assume now that $\M$ is an irreducible weight $\LL$-module such that
there exists $\lambda\in\C$ satisfying ${\rm
dim\,}\M_\lambda=\infty$.
\begin{lemm}
\label{lemm2.3} There exists at most one $i\in\Z$ such that ${\rm
dim\,}\M_{\lambda+i}<\infty$.
\end{lemm}
{\it Proof.~}~Assume that ${\rm dim\,}\M_{\lambda+i}<\infty$ and
${\rm dim\,}\M_{\lambda+j}<\infty$ for different $i,j\in\Z$. Without
loss of generality we may assume $i=1$ and $j>1$. Let $\V$ be the
following $\LL$-module:
\begin{eqnarray}\label{Kernel}
\!\!\!\!\!\!\!\!\!\!\!\!\!\!\!\!\!\!\!\!&\!\!\!\!&\mbox{Ker}(L_1:\M_\lambda\rightarrow\M_{\lambda+1})\cap\mbox{Ker}
(Y_1:\M_\lambda\rightarrow \M_{\lambda+1})\cap\mbox{Ker}
(M_1:\M_\lambda\rightarrow
\M_{\lambda+1})\nonumber\\
\!\!\!\!\!\!\!\!\!\!\!\!&\!\!\!\!\!\!\!\!\!\!\!\!&\cap\,\mbox{Ker}(L_j:\M_\lambda\rightarrow
\M_{\lambda+j}) \cap\mbox{Ker}(Y_j:\M_\lambda\rightarrow
\M_{\lambda+j})\cap\,\mbox{Ker}(M_j:\M_\lambda\rightarrow
\M_{\lambda+j}).
\end{eqnarray}
Since
\begin{eqnarray*}
{\rm dim\,}\M_\lambda=\infty,\ \ \  {\rm
dim\,}\M_{\lambda+1}<\infty\ \ \ {\rm and \ }\ {\rm
dim\,}\M_{\lambda+j}<\infty,
\end{eqnarray*}
all cokernels of the maps appearing in (\ref{Kernel}) are
finite-dimensional. Thus
\begin{equation}\label{infty1}{\rm dim\,}\V=\infty.\end{equation}
Since
\begin{eqnarray*}
&&[L_1,L_p]=(p-1)L_{p+1}\neq 0,\,\,[L_1,M_r]=rM_{r+1}\neq 0 \ {\rm
and}\\
&&[L_1,Y_q]=(q-\frac{1}{2})Y_{q+1}\neq 0
\end{eqnarray*}
for $1\neq p\in\Z,\,q\in\Z,\,1\leq r\in\Z$, we get
\begin{eqnarray}\label{SSSS}
&&L_p\V=0,\ \ \ \ \ p=1,j,j+1,j+2,\cdots,\ \mbox{\ \ and}\nonumber\\
&&Y_q\V=M_r\V=0,\ \ q,r=1,2,\cdots.
\end{eqnarray}
If there would exist $0\neq v\in\V$ such that $L_2v=0$, then
$Y_1v=M_1v=L_1v=L_2v=0$ and $\M$ would be a Harish-Chandra module by
Lemma \ref{lemm2.2}, a contradiction. Hence $L_2v\neq 0$ for all
$v\in\V$. In particular, $${\rm dim\,}L_2\V=\infty.$$ Since ${\rm
dim\,}M_{\lambda+1}<\infty$, there exists $0\neq w\in L_2\V$ such
that $Y_{-1}w=M_{-1}w=L_{-1}w=0$ (as in the proof of
(\ref{infty1})). Let $w=L_2u$ for some $u\in\V$. For all $k\geq j$,
using (\ref{SSSS}), we have
$$L_kw=L_kL_2u=L_2L_ku+(2-k)L_{k+2}u=0.$$ Hence $L_kw=0$ for
all $k=1,j,j+1,j+2,\cdots.$ Since
\begin{eqnarray*}
&&[L_{-1},L_l]=(l+1)L_{l-1}\neq 0,\ [L_l,M_{-1}]=M_{l-1}\neq 0\
\mbox{\,and\,}\\
&&[Y_{-1},L_l]=(1+\frac{l}{2})Y_{l-1}\neq 0
\end{eqnarray*}
for all $l>1$, we get inductively $L_kw=Y_kw=M_kw=0$ for all
$k=1,2,\cdots.$ Hence $\M$ is a Harish-Chandra module by Lemma
\ref{lemm2.2}, again a contradiction. The lemma follows.$\hfill\Box$

Because of Lemma \ref{lemm2.3}, we can now fix the following
notation: $\M$ is an irreducible weight $\LL$-module, $\mu\in\C$ is
such that ${\rm dim\,}\M_\mu<\infty$ and ${\rm
dim\,}\M_{\mu+i}=\infty$ for all $i\in\Z\setminus\{0\}$.
\begin{lemm}
\label{lemm2.4} Let $0\neq v\in\M_{\mu-1}$ be such that
$Y_1v=M_1v=L_1v=0$. Then
\begin{eqnarray*}
&&(1)\ \  (L_1^3-6L_2L_1+6L_3)L_2v=0,\\
&&(2)\ \  Y_kv=M_kv=0,\\
&&(3)\ \  Y_kL_2v=M_kL_2v=0,
\end{eqnarray*}
where $k=1,2,\cdots$.
\end{lemm}

\noindent{\it Proof.~}~(1) is the conclusion of Lemma 4 in
[\ref{MZ}].

(2) Since $Y_1v=M_1v=L_1v=0$ and
$[L_1,Y_l]=(l-\frac{1}{2})Y_{l+1}\neq 0,\,[L_1,M_l]=lM_{l+1}\neq 0$
for $l\geq 1$, we inductively get (2).

(3) follows from
$Y_kL_2=L_2Y_k-(k-1)Y_{k+2},\,M_kL_2=L_2M_k-kM_{k+2}$ and
(2).$\hfill\Box$\vskip7pt

\noindent{\it Proof of Theorem \ref{theo1.3}.~}~ Denote
\begin{eqnarray*}
\V:=\mbox{Ker}\{L_1:\M_{\mu-1}\rightarrow\M_\mu\}\cap
\mbox{Ker}\{Y_1:\M_{\mu-1}\rightarrow\M_\mu\}\cap
\mbox{Ker}\{M_1:\M_{\mu-1}\rightarrow\M_\mu\}.
\end{eqnarray*}
Since ${\rm dim\,}\M_{\mu-1}=\infty$ and ${\rm dim\,}\M_\mu<\infty$,
similar with the discussion of (\ref{infty1}), we have ${\rm
dim\,}\V=\infty$. For any $0\ne v\in\V$, consider the element
$L_2v$. By Lemma \ref{lemm2.2}, $L_2v=0$ would imply that $\M$ is a
Harish-Chandra module, a contradiction. Hence $L_2v\neq 0$, in
particular, ${\rm dim\,}L_2\V=\infty$. This implies that there
exists $w\in L_2\V$ such that $w\neq 0$ and
$Y_{-1}w=M_{-1}w=L_{-1}w=0$. Using
\begin{eqnarray*}
&&[L_{-1},M_{-k+1}]=(1-k)M_{-k}\ \ \ \  {\rm and}\\
&&[L_{-1},Y_{-k+1}]=(\frac{3}{2}-k)Y_{-k},
\end{eqnarray*}
we obtain
$$Y_{-k}w=M_{-k}w=0,\ k=1,2,\cdots.$$
From Lemma \ref{lemm2.4}, we have $(L_1^3-6L_2L_1+6L_3)w=0,$ and
$Y_kw=M_{k}w=0,\ k=1,2,\cdots.$ Therefore,
\begin{eqnarray*}
&&Y_0w=\frac{2}{3}[L_{-1},Y_1]w=0\ \ {\rm and}\\
&&M_0w=[L_{-1},M_1]w=0.
\end{eqnarray*}
The above relations mean that $Y_k,M_k$ act trivially on $\M$ for
all $k\in\Z$. Therefore, $\M$ is simply a module over the Virasoro
algebra $\V ir$. Thus, Theorem \ref{theo1.3} follows from Theorem
\ref{theo2.1}. $\hfill\Box$\vskip7pt

\noindent{\it Proof of Corollary 1.4.~}~Assume that $\M$ is not a
Harish-Chandra module. Then there should exist $i\in\Z$ such that
${\rm dim\,}\M_{\lambda+i}=\infty$. In this case, Theorem 1.3
implies ${\rm dim\,}\M_\lambda=\infty$, a contradiction. Hence $\M$
is a Harish-Chandra module, and the rest of the statement follows
from [\ref{LZ}].$\hfill\Box$

\end{document}